\documentclass[a4paper,10pt]{article}

\usepackage{latexsym}
\usepackage{amsmath}
\usepackage{amssymb}
\usepackage[english]{babel}
\addtocounter{section}{-1}

\usepackage[T1]{fontenc}
\usepackage[utf8]{inputenc}
\usepackage{longtable}
\usepackage{theorem}
\newtheorem{lem}{Lemma}
\newtheorem{prop}{Proposition}
\newtheorem{theo}{Theorem}
\newtheorem{coro}{Corollary}

\def\BP{\noindent  {\it Proof.} }

\def\EP{\hspace*{\fill}$\Box$

\vspace{1ex}
}

\def\BT{\begin{theorem}}
\def\ET{\end{theorem}}
\def\BL{\begin{lemma}}
\def\EL{\end{lemma}}
\def\BCO{\begin{corollary}}
\def\ECO{\end{corollary}}
\def\BX{\begin{example}}
\def\EX{\end{example}}
\def\BIT{\begin{itemize}}
\def\EIT{\end{itemize}}
\def\IT{\vspace{-1ex} \item}

\def\aa{{\rm (a) }}
\def\bb{{\rm (b) }}
\def\cc{{\rm (\hspace{0.1ex}c\hspace{0.1ex}) }}

\def\O{{\cal O}}

\def\S{{\cal S}}
\def\T{{\cal T}}

\begin{document}

\title{{\bf Tychonoff-like Product Theorems\\ for Local Topological Properties}}
\author{\sc Simon Brandhorst {\rm and} Marcel Ern\'e}
\date{\small Faculty for Mathematics and Physics,\\
Leibniz University Hannover, Germany\\ 
{\em e-mail: sbrandhorst@web.de, \
 erne@math.uni-hannover.de}\\[3pt]
}
\maketitle

\begin{abstract}
\noindent  We consider classes $\T$ of topological spaces (referred to as $\T$-spaces) that are stable under continuous images and frequently under arbitrary products.
A local $\T$-space has for each point a neighborhood base consisting of subsets that are $\T$-spaces in the induced topology. A general necessary and sufficient criterion for a 
product of topological spaces to be a local $\T$-space in terms of conditions on the factors enables one to establish a broad variety of theorems saying that a product of spaces has a certain 
local property (like local compactness, local sequential compactness, local $\sigma$-compactness, local connectedness etc.) 
if and only if each factor has that local property, almost all have the corresponding global property, and not too many factors fail a suitable additional condition.
Many of the results admit a point-free formulation; a look at sum decompositions into components of spaces with local properties 
yields product decompositions into indecomposable factors for certain classes of frames like completely distributive lattices or 
hypercontinuous frames.

\vspace{1ex}

\noindent{\bf Mathematics Subject Classification:}\\
Primary: 54B10. Secondary: 06B35, 54D05, 54D30, 54D45.\\[2mm]
\noindent {\bf Key words:}
(locally) compact, (locally) connected, (local) $\T$-space, product space.

\end{abstract}


\section{Introduction}

Undoubtedly one of the most powerful theorems in topology is Tychonoff's theorem, saying that arbitrary products of compact (not necessarily Hausdorff) spaces are again compact.
As observed by Kelley \cite{Kelley} already in the early fifties of the last century, this theorem is equivalent to the full Axiom of Choice (AC), while its restriction to 
Hausdorff spaces or even to the much larger class of sober spaces was later shown to be equivalent to the weaker Prime Ideal Theorem (PIT) (see \L o\'s and Ryll-Nardzewski \cite{LosRyll}, Johnstone \cite{Johnstone}). 
In the present paper we are dealing primarily with arbitrary products of spaces and make permanent use of the surjectivity of the projections; thus, it will be unavoidable to invoke the full strength of (AC),
and we shall do so without particular emphasis. However, if the Continuum Hypothesis (CH) or even the General Continuum Hypothesis (GCH) is involved, this will be mentioned explicitly.

Of course, many variants of compactness have been studied in the past with respect to their stability under the formation of products. 
Since some weaker forms of compactness like paracompactness or countable compactness behave
badly already under the formation of products of two factors (cf. \cite{Dugundji, Novak,counterexamples}), we shall exclude such properties from our present study.
On the other hand, it will be reasonable to include several product-stable properties that are not typically modified compactness properties. A list of relevant topological properties is added at the end of this note.

Our primary concern is the study of {\em local} topological properties and their behavior under the formation of arbitrary products. 
There is a weak and a strong form of localization, which sometimes coincide, but sometimes differ essentially.  
Let $\T$ be a class of topological spaces, referred to as {\em $\T$-spaces}. 
Then, by a \textit{basic $\T$-space}\index{Basic T-space@Basic $\T$-space} we mean a topological space having an open base consisting of subsets that are $\T$-spaces with respect to the induced topology,
and by a \textit{local $\T$-space}\index{Local T-space@Local $\T$-space} a topological space in which every point has a neighborhood base of subsets being $\T$-spaces with respect to their subspace topology.
Of course, every $\T$-basic space is a local $\T$-space, but the converse fails, for example, if $\T$ is the class of compact spaces, whereas both notions agree in case $\T$ is the class of connected spaces.
If we speak of local compactness or local connectedness, respectively, we refer to the above definition for the class $\T$ of compact or connected spaces, respectively. This definition is the one commonly adopted for local compactness in the absence of higher separation axioms (cf.\ \cite{Erne2009, Gierz, Johnstone}), because to require only at least one compact neighborhood for each point would be too weak in order to derive substantial conclusions. However, in the Hausdorff setting of $T_2$-spaces, both notions of local compactness are equivalent. 
Given an infinite cardinal $\kappa$, a union of fewer than $\kappa$ $\T$-spaces will be referred to as a {\em $\T_k$-space}. Our main results are:

\vspace{1ex}

{\em Let $\S,\,\T$ be classes of topological spaces such that $\T$ is closed under continuous images and 
a product of topological spaces is a $\T$-space iff all factors are $\T$-spaces and fewer than $\kappa$ are not $\S$-spaces.
Then a product of topological spaces is a local (basic) $\T$-space iff all factors are local (basic) $\T$-spaces, all but finitely many are $\T$-spaces, and fewer than $\kappa$ are not $\S$-spaces.}

\vspace{1ex}

{\em Assume the class $\T$ is closed under continuous images and arbitrary products.
If a product of spaces is a $\T_{\kappa}$-space then all factors are $\T_{\kappa}$-spaces and for some $\lambda < \kappa$, fewer than $\lambda$ factors are not $\T$-spaces.
The converse  holds if {\rm (GCH)} is assumed and $\kappa$ is a regular limit cardinal or the successor of a regular cardinal.}

\vspace{1ex}

Applying these two general results to various specific classes of topological spaces, we immediately arrive at several known and some new Tychonoff-like theorems concerning global or local properties of topological spaces.

Often the properties under consideration admit a point-free description, that is, a characterization by properties of the open set frames that are \mbox{invariant} under lattice isomorphisms.
This aspect together with the familiar observation that sums of spaces turn into products of their open set frames enables us to establish product decomposition theorems for certain classes of frames, like supercontinuous (=\,completely distributive) lattices and hypercontinuous frames.

\vspace{1ex}

Parts of our results are contained in the first author's 2012 Bachelor thesis.  

\newpage

\section{Set- and order-theoretical preliminaries}
\label{set}

For later use, let us start by recalling a few set-theoretical conventions. 

Each ordinal (number) is regarded as the set of all smaller ordinals. Thus, $\lambda < \kappa$ and $\lambda \in \kappa$ are equivalent statements about ordinals $\kappa$ and $\lambda$. 
We denote by $|A|$ the cardinality of a set $A$; it is the minimal ordinal equipollent to $A$. A cardinal (number) is such a minimal ordinal, equal to its own cardinality. 

Given a a cardinal number $\kappa$, a set $A$ is said to be \textit{$\kappa$-small}\index{kappa-small@$\kappa$-small} if $|A| < \kappa$, 
and, on the contrary, \textit{$\kappa$-large}\index{kappa-large@$\kappa$-large} if $|A| \geq \kappa$. The cardinal successor of $\kappa$ is denoted by $\kappa^+$. 
If $P$ is any property and $(X_i : i\in I)$ is a family of sets or spaces, we say (by slight abuse of language) that {\em fewer than $\kappa$ of the $X_i$'s have property P} if we mean that the set
$\{ i\in I : X_i \mbox{ has } P\}$ is $\kappa$-small (whereas the requirement that the set $\{ X_i : X_i \mbox{ has } P\}$ be $\kappa$-small may be weaker if some or all of the $X_i$'s coincide).
By definition of $\omega$, the least infinite cardinal, ``$\omega$-small'' means ``finite'', and ```fewer than $\omega$'' means ``finitely many''. For the cardinal successor $\omega_1 = \omega^+$ 
(the least uncountable cardinal), ``$\omega^+$-small'' means ``countable''.

A subset $Q$ of a preordered set $P$ is called {\em cofinal} if for each $x\in P$ there exists a $y \in Q$ with $x\leq y$.
The cofinality ${\rm cf}(\kappa)$ of a cardinal number $\kappa$ is defined as the minimal cardinality of a cofinal subset, and $\kappa$ is said to be {\em regular} iff $\kappa= {\rm cf}(\kappa)$. 
In that case, unions of fewer than $\kappa$ sets that are $\kappa$-small are $\kappa$-small.
We write (GCH) to indicate that the Generalized Continuum Hypothesis is assumed, postulating the equation $\kappa^+\! = 2^\kappa$ for all infinite cardinals $\kappa$. 
The special case of the classical Continuum Hypothesis $\omega^+\! = 2^\omega$ is indicated by (CH).

Notice the following properties of infinite sets and cardinals:

\begin{lem}
\label{infinite}
{\rm (1)} For each infinite set $I$ there is a partition into disjoint subsets \mbox{$J_i \subseteq I$} with $|I|=|J_i|$ for all $i \in I$.

{\rm (2)} Each infinite cardinal $\kappa$ is (isomorphic to) the ordinal sum of all smaller ordinals.

{\rm (3)} {\rm (GCH)} assures $\kappa^\lambda = \kappa$ for all infinite cardinal numbers $\kappa, \lambda$ with $\lambda < {\rm cf}(\kappa)$.
\end{lem}

\BP
(1) Note $I\times I = \Sigma_{i \in I}I = \bigcup \{ I \times \{i\}: i\in I\}$ and $|I|=|I\times I|=|\Sigma_{i \in I}I|$; 
so there is a bijection $f:\Sigma_{i \in I}I \rightarrow I$, and one may put $J_i = f[I\times \{ i\}]$.

(2) For $\omega\! \leq \! \nu \! < \! \kappa$, the ordinal sum $\sum_{\iota < \nu} \iota$ is (isomorphic to) a unique ordinal $\lambda_{\nu} < \kappa$
(as $|\sum_{\iota < \nu} \iota| \leq |\nu |^2 = |\nu | < \kappa$),
whence \mbox{$\sum_{\iota < \kappa }\iota = \sup_{\omega\leq \nu < \kappa} \sum_{\iota < \nu} \iota \leq \kappa$.}\\
The reverse inequality is obvious.

For (3), see e.g.\ Jech \cite{Jech}.
\EP

Let us fix a few order-theoretical standard concepts and notations we shall need in due course.
Given a preordered set $P$, a subset $A$ of $P$ and an element $x \in P$,
\begin{itemize} \itemsep0pt
\item[] $\uparrow\! A:=\{y\! \in\! P : \exists\, x\! \in\!  A\, ( x \leq y) \}$ is the \textit{upper set}\index{Upper set} ({\em upset}) generated by\,$A$,
\item[] $\downarrow\! A:=\{y\! \in\! P : \exists\, x\! \in\! X\, ( y \leq x) \}$ is the \textit{lower set}\index{Lower set} ({\em downset}) generated by\,$A$,
\item[] $\downarrow\! x:= {\downarrow\! \{x\}}$ is the \textit{principal filter}\index{Principal filter} generated by $x$,
\item[] $\uparrow\! x:= {\uparrow\! \{x\}}$ is the \textit{principal ideal}\index{Principal ideal} generated by $x$.
\end{itemize}

The \textit{specialization order}\index{Specialization order} on a space $X$ with topology $\O (X)$ is defined by 
\[x \leq y \ \Leftrightarrow  \ \forall \,U \in \O(X)\ ( x\in U \Rightarrow y \in U) \ \Leftrightarrow \ x\in cl\{ y\}.\]
This relation is always a {\em preorder} or {\em quasiorder} (reflexive and transitive). It is an order (antisymmetric) iff $X$ is $T_0$; and it is the identity relation iff $X$ is $T_1$.

Unless otherwise specified, all order-theoretical statements about topological spaces will refer to the specialization order.  
Thus, in a topological space, principal ideals and point closures coincide: $\downarrow \! y =\{x \in X : x \leq y\}= cl\{y\}.$ 
Hence, every open set is an upper set and every closed set is a lower set. 

For a subset $A$ of a space $X$, the {\em saturation} ${\uparrow\! A}$ is the intersection of all neighborhoods of $A$, while the closure $cl\,A$ may properly contain in ${\downarrow\! A}$;
however, in {\em A-spaces}, where arbitrary intersections of open sets are open, $cl\,A = {\downarrow\! A}$. 

It is also helpful to notice that continuous maps $f$ between spaces are monotone with respect to the specialization order, i.e., $x\leq y$ implies $f(x) \leq f(y)$.
The specialization functor, sending a space to the underlying set endowed with the specialization order, preserves initial structures, and in particular products:  

\begin{prop}
\label{produktspezord}
The product of the specialization orders of a family of topological spaces is the specialization order of the product.
\end{prop}

Indeed, the equivalence $\ x \leq y \ \Leftrightarrow  \ \forall\, i \in I \, (x_i \leq y_i)\ $
is just a reformulation of the well known identity \ \mbox{$cl \{ y \}= cl \prod_{i\in I} \{y_i\} = \prod_{i\in I} cl \{ y_i \}$.}

\vspace{1ex}

Given a cardinal number $\kappa$, a preordered set or a topological space is said to be {\em $\kappa$-filtered} or {\em $\kappa$-down-directed} if every $\kappa$-small subset has a lower bound. 
Topologically speaking, a space is $\kappa$-filtered if and only if each $\kappa$-small set of point closures has nonempty intersection. 
The $\omega$-filtered preordered sets are just the down-directed ones. The $\omega$-filtered spaces are the {\em ultraconnected} ones,  
i.e. those nonempty spaces in which the intersection of any two nonempty closed sets is nonempty (see e.g.\ Steen and Seebach \cite{counterexamples}). 
More generally, the $\kappa$-filtered spaces are those in which every $\kappa$-small set of nonempty closed subsets has nonempty intersection ({\em $\kappa$-ultraconnected spaces}). 

\begin{theo}
\label{kfilt}
{\rm (1)} Monotone maps send $\kappa$-filtered sets to $\kappa$-filtered sets.

{\rm (2)} A product of preordered sets is $\kappa$-filtered iff each factor is $\kappa$-filtered.
\end{theo}

\BP
(1) is straightforward: monotone maps preserve lower bounds. 

(2) Suppose $(X_i\! : i\in I)$ is a family of $\kappa$-filtered preordered sets, and \mbox{$X = \prod_{i\in I} X_i$} is their product.
For a $\kappa$-small $A\subseteq X$, each projection set $\pi_i[A]$ is a $\kappa$-small subset of $X_i$, so there exist lower bounds $x_i$ of $\pi_i[A]$, and then $x =(x_i : i\in I)$ is a lower bound of $A$ in $X$.
The other implication is clear by (1).
\EP

\begin{coro}
\label{filtspaces}
{\rm (1)} Continuous maps send $\kappa$-filtered spaces to $\kappa$-filtered spaces.

{\rm (2)} A product of topological spaces is $\kappa$-filtered iff each factor is $\kappa$-filtered.
\end{coro}

Note that a subset of a space is {\em supercompact} (in the sense that it has a dense point) iff it is $\kappa$-filtered for all $\kappa$. Spaces with a base of supercompact open sets are also called {\em B-spaces}, 
and locally supercompact spaces (in which each point has a neighborhood base consisting of supercompact sets) are also termed {\em C-spaces}. 
A-, B- and C-spaces play a central role in the interplay between order and topology (see e.g.\ \cite{Erne1991, Erne2005, Erne2009}). 

\newpage

\section{Products of spaces with local properties}
\label{Prod}

In order to avoid undesired (but trivial) exceptions, all products considered in the sequel are tacitly assumed to be nonempty, and consequently, their factors are nonempty, too. By (AC), all projections
\[\pi_j : \prod_{i\in I} X_i \rightarrow X_j,\ x \mapsto x_j\]
are surjective, and so are all {\em generalized projections}
\[\pi_J :  \prod_{i\in I} X_i \rightarrow \prod_{j\in J} X_j, \ x \mapsto x|_J \ \ (J\subseteq I).\]
For the sake of later use, we note a well-known fact, which holds in arbitrary categories having products; for the special case of topological spaces, see e.g.\ Dugundji \cite{Dugundji}.

\begin{lem}\label{productassociative}
Let $(X_i : i \in I)$ be a family of topological spaces and $\{ J_k : k \in K \}$ a partition of $I$ into nonempty subsets.
Then 
$$\prod_{i \in I} X_i \cong \prod_{k\in K}(\prod_{j \in J_k} X_j).$$
\end{lem}

\vspace{-3ex}

Henceforth, 
\BIT
\IT[] {\em $\T$ always denotes a class of topological spaces (so-called $\T$-spaces) that is closed under images of continuous surjections.}
\EIT 
In particular, the continuity and surjectivity of the generalized projections guarantees that 
\BIT
\IT[] {\em if a product is in $\T$ then so is any subproduct and each factor.} 
\EIT

Recall that a {\em basic $\T$-space} has an open base consisting of $\T$-subspaces, while in a {\em local $\T$-space} every point has a neighborhood base consisting of $\T$-subspaces.

The following easy but fundamental result is due to R.-E.\ Hoffmann \cite{REHoffmann}.

\begin{prop}\label{localTspace}
The image of any local $\T$-space under a continuous open surjection is a local $\T$-space. 
If a product of topological spaces is a local $\T$-space, all factors are local $\T$-spaces and all but finitely many are $\T$-spaces.
The converse holds as well if $\T$ is closed under products.

The same statements are valid with ``local'' substituted by ``basic''.
\end{prop}

By the classical Tychonoff theorem, a product of topological spaces is compact iff all factors are compact. Hence, Proposition \ref{localTspace} immediately provides the well-known fact that
a product is locally compact iff all factors are locally compact and all but finitely many are compact, and similar conclusions for (path-)connectedness instead of compactness. But there are also less familiar
applications of that very flexible proposition. 

For example, calling a a space {\em compactly based} if it has a base of compact open sets, we see that a product of topological spaces is compactly based 
iff each factor is compactly based and only a finite number of the factors are not compact. A {\em spectral space} is compactly based, sober (see e.g.\ \cite{Gierz} or \cite{Johnstone})
and coherent, where coherence means that finite intersections of compact saturated subsets are compact (in particular, the entire space must be compact). 
Up to isomorphism, the open set frames of such spaces are exactly the {\em coherent} or {\em arithmetic frames}, or equivalently,  
the ideal lattices of bounded distributive lattices. Hence, the spectral spaces are the duals of bounded distributive lattices in the classical Stone duality \cite{Sto2}.
A local variant is the following: according to \cite{Gierz}, a {\em stably compact} space is a locally compact, sober and coherent space. Via the patch functor, these spaces are in one-to-one
correspondence with the {\em compact pospaces} (i.e.\ compact topological spaces equipped with an additional closed order), 
and via the open set functor, they are dual to the {\em stably continuous frames} (see \cite{Gierz} for details).
Since a product of spaces is sober or coherent, respectively, iff each factor has the corresponding property, 
we have the following interesting consequences of Proposition \ref{localTspace}:

\begin{coro}
\label{spectral}
\BIT
\IT[{\rm (1)}] A product of spaces is spectral iff all factors are spectral.
\IT[{\rm (2)}] A product of spaces is stably compact iff all factors are stably compact.
\IT[{\rm (3)}] A product of ordered topological spaces is a compact pospace iff all factors are compact pospaces.
\EIT
\end{coro}

Our first major theorem extends Proposition \ref{localTspace} and will serve as the clue for many Tychonoff-like product theorems involving local properties.

\begin{theo}
\label{localSTspace}
Let $\S$ be a class of topological spaces and $\kappa$ an infinite cardinal such that 
a product of topological spaces is a $\T$-space iff all factors are $\T$-spaces and fewer than $\kappa$ of them are not $\S$-spaces.

Then a product of topological spaces is a local $\T$-space iff all factors are local $\T$-spaces, all but finitely many are $\T$-spaces, and fewer than $\kappa$ are not $\S$-spaces.

The same conclusion holds with ``basic'' instead of ``local''.
\end{theo} 

\BP
Let $X=\prod_{i \in I}X_i$ be a local $\T$-space. Then, by Proposition \ref{localTspace}, all factors are local 
$\T$-spaces and all but finitely many are $\T$-spaces. Suppose that for a set $K \subseteq I$ having cardinality $\kappa$, no $X_k$ with $k\in K$ is an $\S$-space. 
Since $\kappa$ is infinite, there exists a partition of $K$ into $\kappa$-large subsets $K_l$, $l\in L$, $|L|=\kappa$ (Lemma \ref{infinite}).
By Lemma \ref{productassociative},  the spaces $Y _l =\prod_{k \in K_l}X_k$ satisfy $\prod_{l\in L} Y_l \cong X$. 
As this is a local $\T$-space, we may again apply Proposition \ref{localTspace}
to see that $Y_l$ is a $\T$-space for all but finitely many $l \in L$.
Since $L$ is infinite, we find indices $l \in L$ such that $Y_l$ is a $\T$-space.
Thus, by the hypothesis on $\S$ and $\T$, the subset $\{k \in K_l : X_k \not\in \S\}$ is $\kappa$-small. 
But $K_l$ is $\kappa$-large, hence some $X_k$ is an $\S$-space -- a contradiction.
 
The other implication is almost routine: suppose each factor $X_i$ is a local $\T$-space, $J = \{i \in I: X_i\not\in \T\}$ is finite and $\{i \in I : X_i \not\in \S\}$ is $\kappa$-small.
Let $U=\prod_{i\in I}U_i$ be a basic open neighborhood of $x\in X\! = \! \prod_{i \in I}X_i$.
Then \mbox{$K =\{i \in I:U_i \neq X_i\}$} is finite. For $i \in J \mathop{\cup} K$, pick a $\T$-neighborhood $V_i \subseteq U_i$ of $x_i$, and set \mbox{$V_i=X_i\in \T$} for $i \in I \setminus (J \mathop{\cup} K)$. 
By the hypotheses on $\S$ and $\T$, the product $\prod_{i\in I}V_i$ is a $\T$-space and is a $\T$-neighborhood of $x$ contained in $U$.

The proof for the ``basic'' case is quite similar.
\EP

To deduce Theorem \ref{localSTspace}, we have used a part of Proposition \ref{localTspace}; of course, the latter in turn is a consequence of the former, by taking simply the special case $\S = \T$. 

Let us note a further consequence of Proposition \ref{localTspace} combined with Corollary \ref{filtspaces}:

{\em A product of topological spaces is locally $\kappa$-filtered iff all factors are locally $\kappa$-filtered and only finitely many are not $\kappa$-filtered.
In particular, a product of topological spaces is locally ultraconnected iff all factors are locally ultraconnected and only finitely many are not ultraconnected.}

\section{Products of unions of $\T$-spaces}
\label{kappaunion}

Let $\kappa$ be an infinite cardinal number. By a {\em $\T_{\kappa}$-space} we mean a topological space that is representable as a union of fewer than $\kappa$ many $\T$-subspaces.
For example, if $\T$ is the class of compact spaces, the $\T_{\omega}$-spaces are still compact, whereas the $\T_{\omega^+}$-spaces are just the $\sigma$-compact ones.

By our general assumption that $\T$ is closed under continuous images, we have:

\begin{lem}
\label{imageTunion}
The image of a $\T_{\kappa}$-space under a continuous map is a $\T_{\kappa}$-space.
\end{lem}

\begin{prop}\label{Tunion1}
If a product of topological spaces is a $\T_{\kappa}$-space then each factor is a $\T_{\kappa}$-space, and for some cardinal $\lambda\! < \!\kappa$, fewer than $\lambda$ factors are not $\T$-spaces.  
\end{prop}
\BP
Consider a product $ X = \prod_{i \in I} X_i$ and suppose $X=\bigcup_{j \in J} Y_j$ for some $\kappa$-small set $J$ and certain $Y_j\in\T$. 
That the factors $X_i$ must be $\T_{\kappa}$-spaces is clear by Lemma \ref{imageTunion}. By way of contraposition, assume that for each $\lambda < \kappa$, at least $\lambda$ many factors are not $\T$-spaces.
Then, in particular, this holds for $\lambda = |J|$. Hence, there exists a subset $K$ of $I$ and a bijection $g : K \rightarrow J$ such that $X_k\not\in \T$ for all $k \in K$. 
For each $k\in K$, the set $Z_k = X_k\setminus \pi_k [Y_{g(k)}]$ must be nonempty, since $Y_{g(k)}$ and so $\pi_k [Y_{g(k)}]$ is a $\T$-space, while $X_k$ is not. 
Putting $Z_i = X_i$ for $i\in I\setminus K$, we conclude that no $x\in \prod_{i\in I} Z_i$ can be contained in any one of the sets $Y_j$ (otherwise $x_k \in \pi_k [Y_j]$ for $k = g^{-1}(j)$), a contradiction.
\EP

\begin{theo}\label{Tunion2}
{\rm (GCH)} Let $\kappa$ be a regular limit cardinal or the successor of some regular infinite cardinal, and assume that the class $\T$ is closed under products.
Then a product of topological spaces is a $\T_{\kappa}$-space iff all factors are $\T_{\kappa}$-spaces and for some cardinal $\lambda < \kappa$, fewer than $\lambda$ factors are not $\T$-spaces.
\end{theo}

\BP
Let $(X_i : i \in I)$ be a family of $\T_{\kappa}$-spaces so that \mbox{$J=\{j \in I : X_j\not\in \T\}$} is $\lambda$-small for an infinite cardinal $\lambda < \kappa$, 
and each $X_j$ is a union of $\T$-subspaces $Y_{k,j}$ ($k \in \lambda_j, \ \lambda_j < \kappa$). 
The supremum of all those $\lambda_j$ is smaller than $\kappa$, by regularity of $\kappa$. Hence, taking $\lambda$ sufficiently large, we may assume $\lambda_j = \lambda$ for all $j \in J$. 
For $f \in \lambda^J$ (the set of all functions from $J$ into $\lambda$) put $Z_f = \prod_{i\in I} W_i$, 
where $W_j = Y_{f(j),j}$ for $j\in J$ and $W_i = X_i$ otherwise. By product closedness of $\T$, each $Z_f$ is a $\T$-space. We have
$\prod_{i \in I} X_i = \bigcup_{f \in \lambda^J} Z_f$, since for each $x=(x_i\! : i\!\in\! I)\in \prod_{i \in I} X_i$ and $j \in J$, there is an index $f(j)$ such that $x_j \in Y_{f(j),j}$, whence $x \in Z_f$. 
It remains to check that $|\lambda^J |< \kappa$ (then $X$ is the union of a $\kappa$-small family of $\T$-subspaces).
Either $\kappa$ is a limit cardinal; in that case, $|\lambda^J| \leq \lambda^\lambda = \lambda^+ < \kappa$ by (GCH); 
or, $\kappa = \rho^+$ for a regular cardinal $\rho$ and then $|\lambda^J| \leq \rho^{|J|} = \rho < \kappa$ ; for $\rho^{|J|} = \rho$, use $|J| < \rho = {\rm cf}(\rho)$ and (GCH) in case $\rho > \omega$.

The converse implication is assured by Proposition \ref{Tunion1}.
\EP

\begin{coro}
\label{Tunion3}
{\rm (GCH)} Let $\kappa$ be the successor of a regular infinite cardinal $\lambda$, and let $\T$ be closed under products. 
Then a product of topological spaces is a local (resp.\ basic) $\T_{\kappa}$-space iff all factors are local (resp.\ basic) $\T_{\kappa}$-spaces, 
all but finitely many factors are $\T_{\kappa}$-spaces, and fewer than $\lambda$ are not $\T$-spaces.
\end{coro}

\BP
Apply Theorem \ref{localSTspace}, with $\T_{\kappa}$ instead of $\T$ and $\T$ instead of $\S$, to the hypothesis assured by Theorem \ref{Tunion2}.
\EP

\section{$\kappa$-union compactness}
\label{Localcomp}

In the following three sections, we have a look at some global and local compactness properties and their productivity, 
in order to demonstrate how Theorems \ref{localSTspace} and \ref{Tunion2} come into play.

Let us call a topological space \textit{$\kappa$-union compact} if it is the union of a $\kappa$-small family of compact subsets. In other words, for the class $\T$ of compact spaces, 
the $\T_{\kappa}$-spaces are just the $\kappa$-union compact ones, and in particular, $\omega^+$-union compactness means $\sigma$-compactness. 
(We avoid the terminology ``$\kappa$-compact'', which in other contexts is reserved to mean that every open cover contains a $\kappa$-small subcover, 
a generalized compactness property we shall not discuss here; for example, ``$\omega^+$-compact'' means ``Lindel\"of''.) 

In this specific setting, Lemma \ref{imageTunion}, Proposition \ref{Tunion1}, Theorem \ref{Tunion2} and Corollary \ref{Tunion3} amount to:

\begin{lem}
\label{imagekappa}
Continuous images of $\kappa$-union compact spaces are $\kappa$-union compact.
\end{lem}

\begin{prop}\label{k-compact2}
If a product of topological spaces is $\kappa$-union compact then all factors are $\kappa$-union compact, and for some $\lambda < \kappa$, fewer than $\lambda$ are not compact. 
\end{prop}

\begin{theo}\label{kappa-union compact}
{\rm (GCH)} Let $\kappa$ be a regular limit cardinal or the successor of some regular infinite cardinal.
Then a product of topological spaces is $\kappa$-union compact iff all factors are $\kappa$-union compact and for a cardinal $\lambda\! <\! \kappa$, fewer than $\lambda$ factors are not compact.
\end{theo}

\begin{coro}
{\rm (GCH)} Let $\kappa$ be the successor of a regular infinite cardinal $\lambda$. 
Then a product of topological spaces is locally $\kappa$-union compact iff all factors are locally $\kappa$-union compact, 
all but finitely many factors are $\kappa$-union compact, and fewer than $\lambda$ are not compact.

In particular, a product of topological spaces is (locally) $\sigma$-compact iff all factors are (locally) $\sigma$-compact and all but finitely many factors are compact.
\end{coro}

An inspection of the proof for Theorem \ref{Tunion2} shows, neither (GCH) nor (CH) is needed for the conclusion about the productivity of (local) $\sigma$-compactness, 
since the equation $\omega^{|J|} = \omega$ holds for $|J| < \omega$, without assuming (CH).

\vspace{1ex}

Entirely analogous results are obtained for ``connected'' or ``path-connected'' instead of ``compact''. Thus, for example, given a regular infinite cardinal $\lambda$, 

{\em a product of topological spaces has at most $\lambda$ (path) components iff all factors have at most $\lambda$ (path) components and all but fewer than $\lambda$ factors are (path) connected.}

\section{$\kappa$-filtered spaces and $\kappa$-sequential compactness}
\label{special}

Recall that a topological space is said to be \textit{sequentially compact}\index{Sequential compactness} if every sequence in the space has a convergent subsequence.
More generally, given an infinite cardinal $\kappa$, we mean by a {\em $\kappa$-sequentially compact space} a topological space 
in which every $\kappa$-sequence $(x_{\nu} : \nu \in \kappa)$ has a convergent $\kappa$-subsequence
\mbox{$(x_{\tau(\nu )} : \nu \in \kappa)$} (where $\tau : \kappa \rightarrow \kappa$ is a strictly monotone increasing map). 

Sequential compactness is incomparable to compactness: while the ordinal space $\omega^+$ is sequentially compact but not compact, 
an uncountable power of the real unit interval is compact but not sequentially compact (cf.\ \cite{counterexamples}).
It is easy to see that the class of $\kappa$-sequentially compact spaces is closed under continuous images.

A preordered set or topological space has been called $\kappa$-filtered if each $\kappa$-small subset has a lower bound. The previous notions are related as follows:

\begin{lem}
\label{kappaconvergence}
$\!$The following conditions on a topological space\,$X$\,are equivalent:
\BIT
\IT[\aa] Every $\kappa$-sequence in $X$ converges.
\IT[\bb] $X$ is $\kappa$-sequentially compact and $\kappa$-filtered.
\IT[\cc] $X$ is $\kappa^+$-filtered.
\EIT
\end{lem}

\BP
(a)$\,\Rightarrow\,$(b)\,and\,(c): Clearly, (a) implies $\kappa$-sequential compactness, so it suffices to verify that (a) yields a lower bound for
any subset $A= \{ a_{\nu} : \nu\in \kappa\}$. By Lemma \ref{infinite}, one has a unique isomorphism $f$ between $\kappa$ 
and the ordinal sum $\sum_{\lambda < \kappa} \lambda = \bigcup \,\{ \lambda \times \{ \lambda\} : \lambda < \kappa\}$ 
(ordered by $(\iota,\lambda)\leq (\iota',\lambda')$ iff $\lambda < \lambda'$ or $\lambda = \lambda'$ and $\iota \leq \iota'$).
Define a $\kappa$-sequence $(x_{\mu} : \mu \leq \kappa)$ by $x_{\mu} = a_{\pi_1 \circ f(\mu )}$, where $\pi_1$ is the projection onto the first coordinate.
This $\kappa$-sequence converges to a point $x$ in $X$. Hence, for any neighborhood $U$ of $x$, 
there is a $\lambda \in \kappa$ such that $U$ contains the set $\{ x_{\mu} : \lambda \leq \mu < \kappa\} = A$. 
(To see that each $a_{\nu}$ is an $x_{\mu}$ with $\mu \geq \lambda$, note that $\kappa$ is a limit ordinal, whence there exists some 
$\lambda'$ with $\max \{ \pi_2 \circ f(\lambda), \nu \} < \lambda' < \kappa$, for which it follows that $f(\lambda ) < (\nu, \lambda' )$, and as $f$ is monotone, $f(\mu) < (\nu, \lambda')$ for all $\mu < \lambda$.
By contraposition, the unique $\mu$ with $f(\mu )= (\nu, \lambda' )$ satisfies $\lambda \leq \mu$ and $x_{\mu} = a_{\pi_1 \circ f(\mu )} = a_{\nu}$.) 
Thus, $x$ is a lower bound of $A$. 

(b)$\,\Rightarrow\,$(c): Let $A = \{ a_{\nu} : \nu \in \kappa\}$ be any $\kappa^+$-small subset of $X$. 
For each (ordinal!) $\iota\in \kappa$, the set $\{ a_{\nu} : \nu \in \iota\}$ has a lower bound $x_{\iota}$, because $X$ is $\kappa$-filtered and $\iota$ ({\it a fortiori} $|\iota|$) is smaller than $\kappa$.
By the hypothesis that $X$ is $\kappa$-sequentially compact, we find a strictly monotone increasing $\tau : \kappa \rightarrow \kappa$ such that the $\kappa$-subsequence
$(x_{\tau(\lambda)} : \lambda \in \kappa )$ converges to a point $x\in X$. Thus, for each open neighborhood $U$ of $x$, there is a $\mu\in \kappa$ such that
$U$ contains each $x_{\tau (\lambda)}$ with $\mu\leq \lambda < \kappa$. Now, for any $\nu\in \kappa$, we find a $\lambda$ with $\mu \leq \lambda < \kappa$ and $\nu < \tau (\lambda )$.
It follows that $x_{\tau (\lambda )}\in U$ and  $x_{\tau (\lambda )}\leq a_{\nu}$, hence $a_{\nu}\in U$ (because $U$ is an upper set). This shows that $x$ is a lower bound of $A$.

(c)$\,\Rightarrow\,$(a): The range of a $\kappa$-sequence $(x_{\mu}\! : \mu\in \kappa)$ in $X$ is a $\kappa^+$-small subset and has therefore a lower bound $x$. 
By definition of the specialization order, every neighborhood of $x$ contains the whole sequence, so it converges to $x$.
\EP

For $\kappa = \omega$, the equivalence of (a) and (b) was also observed by Lipparini \cite{Lipparini}.

\begin{prop}
\label{kappasequ}
If a product of topological spaces is $\kappa$-sequentially compact then fewer than $2^{\kappa}$ (with {\rm (GCH):} at most $\kappa$) factors are not $\kappa^+$-filtered.
\end{prop}

\BP
By way of contraposition, let us look at a product $X = \prod_{i\in I} X_i$ of spaces none of which is $\kappa^+$-filtered, 
where we may assume that $I$ is the set of all functions from $\kappa$ into $\kappa$ (on account of the equation $\kappa^{\kappa} = 2^{\kappa}$).
Using Lemma \ref{kappaconvergence}, pick for each $i\in I$ a non-convergent $\kappa$-sequence $(x^i_\nu : \nu \in \kappa)$ in $X_i$. Consider the elements
$y_{\nu} \in X$ with $y_{\nu,i} = x^i_{i(\nu)}$ (note $i(\nu )\in \kappa$ for $i\in I$). 
Any $\kappa$-subsequence of $y = (y_{\nu}: \nu \in \kappa )$ is of the form $y \circ j = (y_{j(\nu)} : \nu \in \kappa )$ for a strictly monotone increasing $j : \kappa \rightarrow \kappa$.
Since $j$ is injective, there is an $i\in I$ with $i\circ j = id_{\kappa}$. Coordinatewise, this means $y_{j(\nu ), i} = x^i_{i\circ j (\nu )} = x^i_{\nu}$. Thus, no $\kappa$-subsequence of
$y$ can converge, as some of its coordinate sequences do not converge. Hence, $X$ is not $\kappa$-sequentially compact. 
\EP

\begin{theo}
\label{sequ}
{\rm (CH)} {\rm (1)} A product of topological spaces is sequentially compact iff all factors are sequentially compact and all but countably many are $\omega^+$-filtered.

{\rm (2)} A product of topological spaces is locally sequentially compact iff all factors are locally sequentially compact, 
all but finitely many are sequentially compact, and all but countably many are $\omega^+$-filtered.
\end{theo}

\BP
(1) It is well-known that a product of countably many sequentially compact topological spaces is sequentially compact.
Hence, if all other factors are $\omega^+$-filtered, the whole product is still sequentially compact.

For $\kappa = \omega$, Proposition \ref{kappasequ} requires only (CH) and provides the other implication.

(2) Use (1) and apply Theorem \ref{localSTspace} to the class $\S$ of $\omega^+$-filtered spaces and the class $\T$ of sequentially compact spaces.
\EP

An independent proof of (1) under weaker cardinal assumptions was given recently by Lipparini \cite{LippariniII}. 

\vspace{1ex}

It would be interesting to discover how far Theorem \ref{sequ} may be extended to $\kappa$-sequentially compact spaces. One step in that direction is provided by

\begin{coro}
\label{kseqproduct}
If $\kappa$ has cofinality $\omega$\,then a product of $\kappa$-sequentially compact spaces all but countably many of which are $\kappa^+$-filtered is $\kappa$-sequentially compact.
\end{coro}

\BP
By the product-stability of $\kappa^+$-filteredness and Lemma \ref{kappaconvergence}, it suffices to verify that
for any sequence $(X_n : n\in \omega )$ of $\kappa$-sequentially compact spaces, their product $X$ is $\kappa$-sequentially compact, too.
The proof is similar to the classical case of sequential compactness but slightly more involved. 
From any $\kappa$-sequence $x = (x_{\iota} : \iota \in \kappa )$ in $X$, one may extract successively $\kappa$-subsequences $x\circ \psi_n = (x_{\psi_n(\iota)} : \iota \in \kappa )$ such that
\vspace{-1ex}
\begin{itemize}
\item[(1)] $\varphi_n : \kappa \rightarrow \kappa$ and $\psi_n = \varphi_0 \circ ... \circ \varphi_n$ are strictly monotone increasing\\[-3.5ex]
\item[(2)] $n\leq m$ implies $\psi_n (\iota ) \leq \psi_m (\iota )$ (since each $\varphi_k$ is extensive)\\[-3.5ex]
\item[(3)] the coordinate $\kappa$-subsequence $(x_{\psi_n(\iota ),n}\! : \iota\! \in\! \kappa)$ converges to some $z_n$\,in $X_n$.\\[-3.5ex]
\end{itemize}
Pick a monotone increasing cofinal subsequence $(\alpha_n : n\in \omega )$ of $\kappa$, put
$n_\iota : = \min \{ n\in \omega : \iota \leq \alpha_n \}$ for $\iota \in \kappa$ and define a $\kappa$-subsequence $x\circ \varrho = (x_{\varrho (\iota )} : \iota \in \kappa )$ by
$\varrho (\iota) =\psi_{n_\iota} (\iota)$. This $\varrho$ is actually strictly monotone, since by (1) and (2),
$\iota < \lambda < \kappa \,\Rightarrow\, n_{\iota} \leq n_{\lambda} \,\Rightarrow\, \varrho (\iota ) = \psi_{n_{\iota}}(\iota) \leq \psi_{n_{\lambda}} (\iota) < \psi_{n_{\lambda}} (\lambda ) = \varrho (\lambda )$.

In order to assure that $x\circ \varrho$ converges to $z = (z_n : n\in \omega )$ in $X$, it suffices to check that each coordinate $\kappa$-sequence $(x\circ \varrho )_n = (x_{\varrho (\iota), n}\! : \iota\!\in\! \kappa )$ converges to $z_n$. To that aim, consider a neighborhood $U$ of $z_n$ and a $\mu \in \kappa$ such that $x_{\psi_n (\iota ) , n} \in U$ for all $\iota\in \kappa$ with $\mu \leq \iota$. We may assume $\mu > \alpha_n$ and consequently $n < n_{\mu}$ (otherwise $\mu \leq \alpha_{n_{\mu}} \leq \alpha_n$). 
For fixed $\iota$ and the strictly monotone increasing, hence extensive map $\theta = \varphi_{n+1}\circ ...\circ \varphi_{n_{\iota}}$, we get $\varrho (\iota ) = \psi_n \circ \theta (\iota)$.
Thus, $\iota \leq \theta (\iota )$ and therefore $x_{\varrho (\iota ),n} = x_{\psi_n (\theta (\iota)),n} \in U$ for all $\iota$ with $\mu \leq \iota < \kappa$, as desired.
\EP

\section{Hypercompact and supercompact spaces}
\label{hyper}

By evident reasons, we  call a finitely generated upper set ${\uparrow\! F}$ in a preordered set a {\em foot}. Referring to the specialization order, 
a subset $H$ of a topological space $X$ is called \textit{hypercompact}\index{Hypercompact} if its saturation is a foot, i.e., ${\uparrow \! H} = {\uparrow \! F}$ for some finite $F$ (contained in $H$);
if $F$ can be chosen to be a singleton, $H$ is called \textit{supercompact}\index{Supercompact} (see \cite{Erne1991, Erne2005, Erne2009}).

As demonstrated in \cite{Erne2009}, these notions admit point-free characterizations. 
An element $c$ of a poset $T$ is said to be {\em hypercompact} (resp.\ {\em supercompact} or {\em completely join-prime}) if $T\setminus {\uparrow\!c}$ is a finitely generated lower set (resp.\ a principal ideal). 
Clearly, these properties are stronger than order-theoretical compactness of $c$, which means that $T\setminus {\uparrow\!c}$ is closed under directed joins (cf.\ \cite{Erne2009, Gierz, Johnstone}). 
Now, it turns out that an open subset of a topological space $X$ is hyper- or super\-compact, respectively, 
if and only if it has the synonymous order-theoretical property, regarded as an element of the open set frame $\O (X)$. 

More to the point for us, not only hyper- and supercompactness, but also local hyper- and supercompactness admit elegant point-free descriptions -- in contrast to local compactness, which is {\em not} invariant under
lattice isomorphisms between the open set frames (see \cite{Gierz} for a sophisticated counterexample):
the supercontinuous frames (=\,completely distributive lattices) are, up to isomorphism, the open set frames of locally supercompact spaces \cite{Erne1991, Erne2009}, while the
\mbox{hypercontinuous} frames are the open set frames of locally hypercompact spaces (also called {\em finitely bottomed spaces}; cf.\ \cite{Erne2009, Lawson1985}).
There are several equivalent definitions of hypercontinuity. The most topologically inspired one is this (see \cite{Gierz}): denoting by $\upsilon P$ the {\em upper} or {\em weak topology} 
generated by the complements of the principal ideals of a lattice or preordered set $P$, \textit{hypercontinuity} of $P$ 
means that for each $y \in P$,  the set $\{x \in P: y \in int_{\upsilon P}\, {\uparrow\! x}\}$ is directed and has the join $y$; in the case of a complete lattice $P$, the directedness condition is automatically fulfilled.
A straightforward verification confirms:

\begin{lem}
\label{hyperbildabg}
The class of hypercompact spaces is closed under the formation of continuous images.
\end{lem}

\begin{prop}
\label{foot}
A product of preordered sets is a foot iff all factors are feet and all but finitely many have least elements. 
\end{prop}

\BP
Let $\prod_{i\in I}P_i$ be a foot and find a finite $F\!\subseteq\! \prod_{i\in I}P_i$ with ${\uparrow \! F} = \prod_{i\in I}P_i$. Then each factor $\uparrow \! \! \pi_i(F)=P_i$ is a foot. 
Assume that infinitely many factors do not have least elements. Thus, if \mbox{$F=\{f_1, ... , f_n\}$}, there are distinct $i_1, ... i_n \in I$ such that $P_{i_k}$ has no least element. 
This implies that $Y_k = P_{i_k} \setminus \! \! \uparrow \! \! \pi_{i_k}(f_k)$ is not empty. Put $Y_i = X_i$ for $i\in I \setminus \{ i_1,...,i_n\}$. 
Then no $x\in \prod_{i\in I} Y_i$ can be contained in $\uparrow \! F$, which is a contradiction.

Conversely let $(P_i : i \in I)$ be a family of feet such that $J =\{i\in I : P_i$ has no least element$\}$ is finite. 
For each $i\in I$ pick $F_i \subseteq P_i$ with minimal cardinality such that $\uparrow \! F_i = P_i$. 
Since almost all $F_i$ are singletons, $|\prod_{i\in I}F_i|=|\prod_{i\in J}F_i|< \infty$, and  \mbox{$\uparrow \prod_{i\in I}F_i = \prod_{i\in I} \uparrow \! \! F_i = \prod_{i\in I}P_i$}. 
We conclude that $\prod_{i\in I}P_i$ is a foot.
\EP

Translating the last result into the language of topological spaces, we conclude: 

\begin{theo}
\label{hyperproduct2}
A product of topological spaces is hypercompact iff all factors are hypercompact and all but finitely many are supercompact.
\end{theo}

Now, a further application of Theorem \ref{localSTspace} yields:

\begin{coro}
\label{hyperproduct3}
A product of topological spaces is locally hypercompact iff all factors are locally hypercompact and all but finitely many are supercompact.
\end{coro}

An analogous conclusion holds for ``supercompact'' instead of ``hypercompact''. Notice that arbitrary products of supercompact spaces are supercompact.

\section{$\kappa$-hypercompact spaces}
\label{kappahyper}

A common generalization of hypercompactness and supercompactness is provided by the following definition (see \cite{Erne2009}). 
Let $\kappa$ be a cardinal number. A subset $H$ of a topological space $X$ is called \textit{$\kappa$-hypercompact}\index{kappa-hypercompact@$\kappa$-hypercompact} iff there exists a $\kappa$-small subset $F \subseteq X$ such that ${\uparrow \! F} = {\uparrow\! H}$. In particular, the whole space is $\kappa$-hypercompact iff it is the saturation of a $\kappa$-small subset. 
By definition, ``$\omega$-hypercompact'' means ``hypercompact'', and ``$2$-hypercompact'' means ``super\-compact''. 
One easily verifies that  Lemma \ref{hyperbildabg} extends to $\kappa$-hypercompact spaces:

\begin{lem}
\label{kappahyperbildabg}
The class of $\kappa$-hypercompact spaces is closed under the formation of continuous images and of closed or at least lower subsets. 
\end{lem}

For the class $\T$ of supercompact spaces, the $\T_{\kappa}$-spaces are just the $\kappa$-hypercompact ones. Hence, a further application of Proposition \ref{Tunion1} and Theorem \ref{Tunion2} gives:

\begin{prop}
\label{khyper-fast}
If a product of spaces is $\kappa$-hypercompact then so is each factor, and for some cardinal $\lambda\! < \!\kappa$, fewer than $\lambda$ factors are not supercompact.
\end{prop}

\begin{theo}
{\rm (GCH)} Let $\kappa$ be a regular limit cardinal or the successor of a regular infinite cardinal.
Then a product of topological spaces is $\kappa$-hypercompact iff all factors are $\kappa$-hyper\-compact and for some $\lambda < \kappa$, fewer than $\lambda$ factors are not supercompact.
\end{theo}

Finally, Theorem \ref{localSTspace} with $\S$ the class of supercompact spaces and $\T$ the class of $\kappa$-hypercompact spaces amounts to:

\begin{coro}
{\rm (GCH)} Let $\kappa$ be the successor of a regular infinite cardinal $\lambda$. A prod\-uct of spaces is locally $\kappa$-hypercompact iff each factor is locally $\kappa$-hypercompact, 
all but finitely many are $\kappa$-hypercompact, and fewer than $\lambda$ are not supercompact.
\end{coro}

\section{Sum decompositions of spaces and product decompositions of frames}
\label{sum} 

One useful effect of point-free thinking is the observation that sum decompositions of spaces into open components yield product decompositions of the corresponding frames into indecomposable factors:
$$\textstyle{\O (\sum_{i\in I} X_i ) \simeq \prod_{i\in I} \O (X_i).}$$
The best known class of topological spaces in which all components are (closed and) open is that of locally connected spaces, 
which alternatively may be described as basic $\T$-spaces or as local $\T$-spaces for the class $\T$ of connected spaces. 
We are now going to single out a few specific classes of such spaces.

\begin{lem}
\label{hyperfinitecomponents}
A hypercompact topological space has only finitely many components, hence a unique finite sum decomposition into indecomposable summands.
\end{lem}

\BP
Since components are closed, two points in different components have disjoint closures, hence no common lower bound in the specialization order. 
Thus, the union of an infinite number of components cannot be a foot.
\EP

As in the case of locally compact spaces, but by different arguments, one has:

\begin{prop}
\label{openclosed}
Any intersection of a closed or at least lower set and an open set of a locally {\rm ($\kappa$-)}hypercompact space is locally {\rm ($\kappa$-)}hypercompact. 
\end{prop}

\BP
Let $U$ be an open and $A$ a lower set in a locally ($\kappa$-)hypercompact space $X$. For $x\! \in \! U \cap A$ and a neighborhood $V$ of $x$ in $X$, we find a ($\kappa$-)hypercompact neighborhood $H \subseteq V \cap U$ of $x$. Now, $H \cap A$ has the desired properties as it is a ($\kappa$-)hypercompact neighborhood of $x$ in $U \cap A$ contained in $U \cap A \cap V$.  
\EP

\begin{theo}
\label{locallyhyperconnected}
Each point of a locally hypercompact topological space has a neighborhood base of connected hypercompact sets and is therefore locally connected.
Hence, every locally hypercompact space has a unique sum decomposition into locally hypercompact and sum-indecomposable subspaces (its components).
\end{theo}

\BP
Let $X$ be a locally hypercompact space and  $x \in U\in \O (X)$. Then we find a hypercompact neighborhood $H \subseteq U$ of $x$. By Lemma \ref{hyperfinitecomponents}, 
it has only finitely many connected components, which are therefore clopen. By $G$ we denote the component of $x$ in $H$. 
Since G is a relatively closed subset of a hypercompact set, it is hypercompact. As it is relatively open in $H$, we find a $V \in \O (X)$ with $G=H\cap V$. %
Now $H$ is a neighborhood of $x$, so $x$ is in its interior. We deduce that $x \in \mbox{\it int}\,H \cap V \subseteq G \subseteq H \subseteq U$.
Thus, $G$ is a connected hypercompact neighborhood contained in $U$. 
\EP

As mentioned earlier, the hypercontinuous (resp.\ supercontinuous) frames are isomorphic to locally hypercompact (resp.\ supercompact) topologies.  
Hence, Theorem \ref{locallyhyperconnected} immediately yields:

\begin{coro}
\label{producthypercont}
Every hypercontinuous (respectively, supercontinuous) frame has a product decomposition into product-indecomposable hypercontinuous (respectively, supercontinuous) factors.
\end{coro}

Since locally connected spaces also admit a point-free description, one has a similar product decomposition for frames in which every element is a join of connected ones, where an element $c$ of a bounded lattice is 
{\em connected} iff $a\vee b = c$ and $a \wedge b = 0$ imply $a=c$ or $b=c$.

\section{Table: Productivity of topological properties}
\label{table}

\setlongtables
\begin{longtable}{| c | c | c | c |}
\hline 
 \hspace{2ex} {\em stable under products\,?} \hspace{2ex} & \hspace{2ex} {\em finite} \hspace{2ex} & {\em countable} & {\em arbitrary} \\
\hline
\hline
compact					  & yes	& yes & yes\\
\hline
countably compact	      & no	& no  & no  \\
\hline
paracompact				  & no	& no  & no  \\
\hline
Lindel\"of                & no  & no  & no  \\
\hline
sequentially compact	  & yes	& yes & no  \\
\hline
$\sigma$-compact		  & yes	& no  & no \\
\hline
supercompact			  & yes	& yes & yes\\
\hline
hypercompact			  & yes	& no  & no \\
\hline
connected                 & yes & yes & yes\\
\hline
path-connected            & yes & yes & yes\\
\hline
($\kappa$-)ultraconnected & yes & yes & yes\\
\hline
$T_i$-space ($i=0,1,2,3$) & yes & yes & yes\\
\hline
$T_4$-space (normal)      & no  & no  & no \\
\hline
sober                     & yes & yes & yes\\
\hline
spectral                  & yes & yes & yes\\
\hline
stably compact            & yes & yes & yes\\
\hline
\end{longtable}

\end{document}